\documentstyle{mathnachtmp}
\Year{1997}
\Received{ }           %Please insert the date into the brackets.
\input epsf
\begin{document}
%%%%%%%%%%%%%%%%%%%%%%%%%%%%%%%%%%%%%%%%%%%%%%%%%%%%%%%%%%%%%%%%%%%%%%%%%%%%%%%
\Title{On the Spectrum of the Orr-Sommerfeld Equation on the semiaxis}
\Shorttitle{ The Spectrum of the Orr-Sommerfeld Equation} 
\By{{\sc J.~Lahmann}, {\sc M.~Plum}} 
\Names{J.~Lahmann, M.~Plum}
\Dedicatory{ }
\Subjclass{ } 
\Email{Jan.Lahmann@math.uni-karlsruhe.de}
\maketitle
%%%%%%%%%%%%%%%%%%%%%%%%%%%%%%%%%%%%%%%%%%%%%%%%%%%%%%%%%%%%%%%%%%%%%%%%%%%%%%%
%                                                                             %
%                   Please insert now the article body.                       %
%                                                                             %
%%%%%%%%%%%%%%%%%%%%%%%%%%%%%%%%%%%%%%%%%%%%%%%%%%%%%%%%%%%%%%%%%%%%%%%%%%%%%%%
\begin{abstract}
  The Orr-Sommerfeld equation is a spectral problem which is known to
  play an important role in hydrodynamic stability. For an appropriate
  operator theoretical realization of the equation, we will determine
  the essential spectrum, and calculate an enclosure of the set of all
  eigenvalues by elementary analytical means.
\end{abstract}

\newsection{Introduction}

The {\it Orr-Sommerfeld} equation is an non-selfadjoint general
eigenvalue problem of the form
\begin{equation} \label{1}
(-D^2+a^2)^2u+ iaR[V\cdot (-D^2+a^2)u +V''\cdot u]=\lambda (-D^2+a^2)u
\quad\mbox{on $I$},
\end{equation}
subject to Dirichlet boundary conditions for $u$ on the real interval
$I$.  Here, $D:=d/dx$, $i$ is the imaginary unit, and $R>0$ is the
Reynolds number of an underlying fluid which moves in a stationary
flow (perpendicular to $I$) with given real-valued flow profile $V\in
C^2(I)$.  This flow is perturbed by a single-mode perturbation with
wave number $a>0$, and the physical question of stability or
instability of the underlying flow in response to this perturbation
arises. This question is closely related to the spectrum of (an
appropriate operator-theoretical realization of) the Orr-Sommerfeld
equation (\ref{1}). Essentially the flow is unstable (with respect to
the wave number $a$) if (\ref{1}) has an eigenvalue with negative real
part.

In this article, we will exclusively be concerned with the case
$I=[0,\infty)$ corresponding to the half-plane flow along a wall (an
overview of results for this case can be found in \cite{herron1}). So
the boundary conditions to be added to (\ref{1}) read
\begin{equation} 
\label{2}
u(0)=u'(0)=u(\infty)=u'(\infty)=0.
\end{equation}
A flow profile $V$ which is of particular interest in this case (but
not under exclusive concern here) is the {\it Blasius} profile defined
by $V:=f'$, where $f$ is the solution of the nonlinear boundary value
problem (Blasius equation, a special case of the Falkner--Skan
equation)
\begin{equation} \label{3}
2f'''+ff''=0\quad\mbox{on $[0,\infty)$},\qquad f(0)=f'(0)=0,\lim\limits_
{x\to\infty}f'(x)=1,
\end{equation}
which can be shown to exist and to be unique, and moreover, to provide
(for $V=f'$)
\begin{eqnarray}
&V\to 1,\ V'\to 0,\ V''\to 0\mbox{\ (exponentially) as 
$x\to\infty$},&\nonumber \\[-0.5\baselineskip]
\label{4}\\[-0.5\baselineskip]
&V>0,\ \ V'>0,\ \ V''<0\mbox{\ on $(0,\infty)$}.&\nonumber
\end{eqnarray}

{ \vspace{-2cm}
\begin{figure}[h]
%\begin{center}
  \epsfverbosetrue \hspace{1.0cm} \epsfxsize=10cm \epsfbox{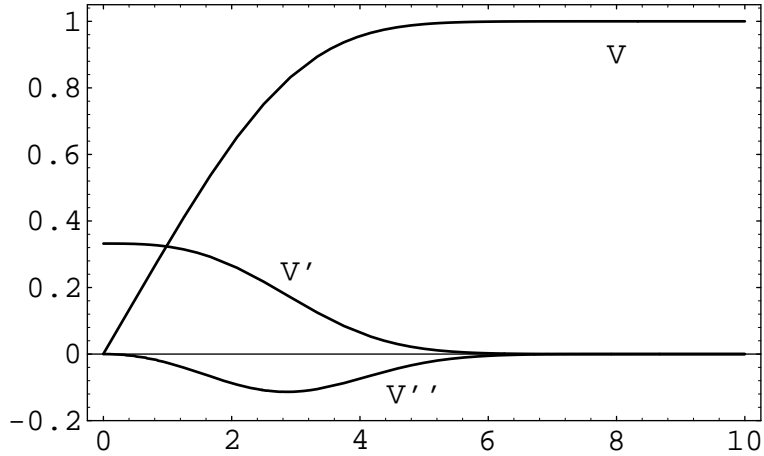}
%\end{center}
  \vspace{-2cm}
\caption{Blasius profile together with it's first two derivatives}
\label{bild0}
\end{figure}
}

We choose the following operator theoretical setting for problem
(\ref{1}), (\ref{2}), which we believe is the most natural and simple
one:

Let $X$ and $Y$ denote the complex Banach spaces $H^2(0,\infty)\cap
H_0^1(0,\infty)$ and $L^2(0,\infty)$, respectively, and let
$D(A):=H^4(0,\infty)\cap H_0^2(0,\infty)$. The operators
$A:D(A)\subset X\to Y$ and $B:X\to Y$ are defined by $Au$ and $Bu$
denoting the left-hand and the right-hand side (without $\lambda$) of
equation (\ref{1}), respectively. It can be shown that $A$ is closed
(see \cite[Theorem IX.9.5]{evans1}, plus the remark that the norm in
$X$ is stronger than the one in $L^2(0,\infty)$), and that $B$ is
one-to-one and onto, and bounded with bounded inverse. Our formulation
of (\ref{1}), (\ref{2}) now reads
\begin{equation} \label{5}
Au\ =\ \lambda Bu.
\end{equation}
So we choose a {\it direct} operator theoretical realization of the
Orr-Sommerfeld problem, rather than a formulation via an auxiliary
operator constructed by Riesz's representation lemma, as done e.g. in
\cite{fischer1}.

Our goal in this article is twofold: In Section 2, we will determine
the essential spectrum of problem (\ref{5}) exactly; the main tool is
the well-known theorem on preservation of the essential spectrum under
relatively compact perturbations. In Section 3, an enclosure for the
set of all eigenvalues of problem (\ref{5}) is calculated by
elementary analytical means. Here, exploiting sign restrictions on the
profile $V$ and its derivatives (such as (\ref{4})), we obtain
results which extend enclosures known in the literature (see
\cite{fischer1}, \cite{miklav1}), even if those have been obtained by
more involved methods. In \cite{joseph1} corresponding results have
been obtained for the case of a compact interval $I$. For further
results on the Orr--Sommerfeld equation on a compact
interval see \cite{shkal1} and the literature cited there.

\newsection{The essential spectrum}

Since several definitions of essential spectra are around in the
literature (see \cite[Chapter IX]{evans1} for an overview), we start
with the definition which we are using.

With $X,Y$ denoting two complex Banach spaces, and $A:D(A)\subset X\to
Y$ a closed linear and $B:X\to Y$ a bounded linear operator, we call
the set
$$
\sigma_{\mbox{ess}}(A,B):=\{\lambda\in\C\, :\, A-\lambda B\mbox{\ 
  is not a Fredholm operator of index $0$}\}
$$
the {\it essential spectrum} of the problem $Au=\lambda Bu$ or of
the pencil $(A,B)$. Here, a closed linear operator $T:D(T)\subset X\to
Y$ is called a {\it Fredholm operator} if its range is closed and has
finite codimension $d(T)$ in $Y$, and if its nullspace has finite
dimension $n(T)$; it has index $0$ if $d(T)=n(T)$. Furthermore,
$$
\varrho(A,B):=\{\lambda\in\C\, :\, A-\lambda B\mbox{\ is one-to-one
  and onto}\}
$$
denotes the {\it resolvent set} of the problem $Au=\lambda Bu$ (or
of the pencil$(A,B)$), and $\sigma(A,B):=\C\setminus\varrho(A,B)$ its
(total) {\it spectrum}. It is quite obvious that
$\sigma_{\mbox{ess}}(A,B)\subset \sigma(A,B)$ and that
$\sigma(A,B)\setminus\sigma_{\mbox{ess}}(A,B)$ consists only of the
eigenvalues of finite geometric multiplicity.

We return to the specific choice for $X,Y,A,B$ introduced in Section
1, i.e., describing the Orr-Sommerfeld problem. Our main result of the
present section is
\begin{theorem}
  Suppose that the flow profile $V$ satisfies $V\to 1, V''\to 0$ as
  $x\to\infty$ $($compare $($\ref{4}$))$. With $R$ and a denoting the
  parameters in $($\ref{1}$)$, the essential spectrum of $($\ref{5}$)$
  is given by
\begin{equation} \label{6}
\sigma_{\mbox{ess}}(A,B)=\{\mu +a^2+iaR\, :\, \mu\in [0,\infty)\}.
\end{equation}
\end{theorem}

\begin{proof}
  We introduce the auxiliary operators
\begin{tabbing}
\quad \= $A_0:D(A)\to Y$, \= $A_0u:=(-D^2+a^2)^2u+ iaR(-D^2+a^2)u,$ \\[0.8 ex]
      \> $K\,:\,D(A)\to Y$,   \> $Ku\,:=\,iaR[(V-1)\cdot 
(-D^2+a^2)u+V''\cdot u]$ \end{tabbing}
so that $A=A_0+K$. Since $V-1$ and $V''$ tend to zero and $K$ is of lower
order than $A_0$, it follows from \cite[Theorem IX.8.2]{evans1} that $K$ 
is relatively
compact to $A_0$, i.e., for each $\|\cdot\|_X$-bounded sequence $(u_n)$ in
$D(A)$ such that $(A_0u_n)$ is bounded in $Y$, $(Ku_n)$ contains a
convergent subsequence. Strictly speaking, the cited Theorem in \cite{evans1}
provides the result with $X$ and $D(A)$ replaced by $L^2(0,\infty)$ and
$H^4(0,\infty)$, respectively, but from this it is easy to obtain the same
statement with $X$ and $D(A)$.

Since $B:X\to Y$ is one-to-one, onto, and bounded with bounded
inverse, the above result implies that $B^{-1}K:D(A)\to X$ is
relatively compact to $B^{-1}A_0: D(A)\to X$. Now the invariance of
the essential spectrum under relatively compact perturbations (see
e.g. \cite[Theorem IX.2.1]{evans1}) provides \begin{equation}
  \label{7}
  \sigma_{\mbox{ess}}(B^{-1}A_0+B^{-1}K,id_X)=\sigma_{\mbox{ess}}
  (B^{-1}A_0,id_X).
\end{equation}
The cited theorem in \cite{evans1} is formulated for densely defined
operators, but 
% restricting $X$ to
% $\tilde{X}:=\overline{D(A)}=H_0^2(0,\infty)$ we obtain (\ref{7}) also
% in our situation where $D(A)$ is not dense in $X$.) 
the proof clearly does not use this assumption. (This density is a 
general assumption in \cite[Chapter IX]{evans1} to be able to work with 
adjoint 
operators, which however do not occur in Theorem IX.2.1 or in its proof.)
Since $B$ and
$B^{-1}$ are bounded, so that $A-\lambda B$ is a Fredholm operator of
index $0$ if and only if $B^{-1}A-\lambda id_X$ is, we obtain
furthermore
$\sigma_{\mbox{ess}}(A,B)=\sigma_{\mbox{ess}}(B^{-1}A,id_X)$, so that
(\ref{7}) and the identity $A_0+K=A$ provide
\begin{equation} \label{8}
\sigma_{\mbox{ess}}(A,B)=\sigma_{\mbox{ess}}(B^{-1}A_0,id_X).
\end{equation}
It remains to be shown that $\sigma_{\mbox{ess}}(B^{-1}A_0,id_X)$
equals the right-hand side of (\ref{6}) which we call $M$ from now on.
This is done in two steps: We prove
\begin{itemize}
\item[a)] $\C\setminus M\subset\varrho(B^{-1}A_0,id_X)$ (implying
  $M\supset\sigma(B^{-1}A_0,id_X)\supset\sigma_{\mbox{ess}}(B^{-1}A_0,id_X)$).
\item[b)] $M\subset\sigma_{\mbox{ess}}(B^{-1}A_0,id_X)$.
\end{itemize}
First we observe that, with $E(x):=e^{-ax}$,
\begin{equation} \label{9}
B^{-1}A_0u = -u''+(a^2+iaR)u+u''(0)\cdot E\quad\mbox{for $u\in D(A)$},
\end{equation}
since the expression on the right-hand side is in $X$ (for $u\in
D(A)$) and $B$ applied to it equals $A_0u$.

{\bf ad a)} Let $\lambda\in\C\setminus M$, so that
$\lambda=\mu+a^2+iaR$ with $\mu\in\C\setminus [0,\infty)$. Thus,
$\mu\in \varrho(C,id_Y)$ where $C:D(C)\subset Y\to Y$ is given by
$D(C):=H^2(0,\infty)\cap H_0^1(0,\infty)$,
$Cu:=-u''$. \\
To prove that $B^{-1}A_0-\lambda id_X: D(A)\to X$ is onto, let $r\in
X$ and define $v:=(C-\mu)^{-1}r\in D(C)=X, w:=(C-\mu)^{-1}E\in
D(C)=X$. It is easy to calculate $w$ in closed form and to show that
$w'(0)\not= 0$.  Therefore, $u:=v-\frac{v'(0)}{w'(0)}w\in X$ solves
the equation
\begin{equation} \label{10}
-u''-\mu u+\frac{v'(0)}{w'(0)}E\ =\ r
\end{equation}
and moreover, $u(0)=u'(0)=0$. Since $u,E,r\in H^2(0,\infty)$,
(\ref{10}) implies $u''\in H^2(0,\infty)$. Therefore $u\in
H^4(0,\infty)$. Altogether, $u\in D(A)$. Since $u(0)=r(0)=0$ and
$E(0)=1$, (\ref{10}) shows that $u''(0)=v'(0)/w'(0)$, and (\ref{10}),
(\ref{9}) then yield $(B^{-1}A_0-
\lambda)u=r$. \\
To prove that $B^{-1}A_0-\lambda id_X$ is one-to-one, let $u\in D(A)$
satisfy $(B^{-1}A_0-\lambda)u=0$, i.e.,
$$
-u''-\mu u+u''(0)E=0.
$$
Thus, $u\in D(C)$ and $(C-\mu)u=-u''(0)E$, i.e., $u=-u''(0)w$.
Since $u'(0)=0, w'(0)\not= 0$, this implies $u''(0)=0$ and therefore
$u=0$.

{\bf ad b)} Let $\lambda\in M$, so that $\lambda=\mu+a^2+iaR$ with
$\mu\in [0,\infty)$. Let $\gamma\in\C$ and $r:=\gamma(1-E)E\in X$, and
suppose that, for some $u\in D(A), (B^{-1}A_0-\lambda)u=r$, so that
(\ref{9}) yields
\begin{equation} \label{11}
-u''-\mu u+u''(0)E = \gamma(1-E)E.
\end{equation}
We will prove that (\ref{11}), together with the fact $u\in D(A)$,
implies $\gamma=0$ and $u=0$, which shows that $B^{-1}A_0-\lambda
id_X$ is one-to-one but not onto. Consequently, it is no Fredholm
operator of index
zero, which implies $\lambda\in \sigma_{\mbox{ess}}(B^{-1}A_0,id_X)$. \\
The general solution of (\ref{11}) satisfies
\begin{equation} \label{12}
u = c_1\varphi_1+c_2\varphi_2+\frac{u''(0)-\gamma}{a^2+\mu}E+\frac{\gamma}
{4a^2+\mu}E^2
\end{equation}
with $c_1,c_2\in\C$, where $\varphi_1(x)=\cos(\sqrt{\mu}x)$ and
$\varphi_2(x)=\sin(\sqrt{\mu}x)$ if $\mu>0,$ $\varphi_1(x)=1$ and
$\varphi_2(x)=x$ if $\mu=0$. The condition $u\in L^2(0,\infty)$
requires $c_1=c_2=0$, and $u(0)=u'(0)=0$ then implies
$$
\frac{\gamma -u''(0)}{a^2+\mu} = \frac{\gamma}{4a^2+\mu},\quad
a\frac{\gamma -u''(0)}{a^2+\mu} = 2a\frac{\gamma}{4a^2+\mu}
$$
which indeed yields $\gamma=0, u''(0)=0$, so that (\ref{12})
provides $u\equiv 0$.
\end{proof}

It is easy to generalize Theorem \ref{6} for flow profiles which
satisfy $V\to c\in\R$ instead of $V\to 1$.
\begin{corollary}
  Suppose that the flow profile $V$ satisfies $V\to c\in\R, V''\to 0$
  as $x\to\infty$. Then the essential spectrum of $($\ref{5}$)$ is
  given by
\begin{equation} \label{6b}
\sigma_{\mbox{ess}}(A,B)=\{\mu +a^2+iaR\,\,c\, :\, \mu\in [0,\infty)\}.
\end{equation}
\end{corollary}

\newsection{An enclosure for the set of all eigenvalues}

In this final section, we will calculate a set enclosing all
eigenvalues of problem (\ref{5}). Since this set will also contain the
straight line which was identified as $\sigma_{\mbox{ess}}(A,B)$ in
Theorem 2.1, it therefore encloses the {\it total} spectrum
$\sigma(A,B)$.

We wish to put emphasize on the simplicity of the methods we use,
which nevertheless provide more accurate enclosures than those known
in the literature \cite{fischer1}, \cite{miklav1}, if sign
restrictions such as (\ref{4}) are exploited.

Let real constants $V_{\min}, V_{\max}, |V'|_{\max}, V''_{\min},
V''_{\max}$ be given such that
\begin{equation} \label{13}
\ \ V_{\min}\le V(x)\le V_{\max}, |V'(x)|\le |V'|_{\max},
V''_{\min}\le V''(x)\le V''_{\max}\quad\mbox{for $x\in [0,\infty)$}
\end{equation}
for the flow profile $V$. Let $\langle\cdot ,\cdot\rangle$ and
$\|\cdot\|$ denote the usual inner product and norm in
$L^2(0,\infty)$.

For any eigenpair $(u,\lambda)\in D(A)\times\C$ of problem (\ref{5})
we obtain (note that $B$ is positive definite)
\begin{equation} \label{14}
\lambda\ =\ \frac{\langle Au,u\rangle}{\langle Bu,u\rangle}
\end{equation}
which is the basis of our enclosures. Partial integration and the
condition $u\in D(A)=H^4(0,\infty)\cap H_0^2(0,\infty)$ provide
\begin{eqnarray*}
\langle Au,u\rangle & = & 
\langle (-D^2+a^2)^2u,u\rangle +iaR\langle V\cdot(-D^2+a^2)u+V''\cdot 
u,u\rangle \\
& = & \|(-D^2+a^2)u\|^2+iaR[\langle Vu',u'\rangle -\langle V'u,u'\rangle 
+a^2\langle Vu,u\rangle ] 
\end{eqnarray*}
so that (\ref{14}) yields
\begin{equation} \label{15}
\lambda\ =\ \beta_1+iaR(\beta_2-\beta_3)
\end{equation}
where
\begin{eqnarray} \label{16}
\beta_1 & := & \frac{\|Bu\|^2}{\langle Bu,u\rangle }\ \ge\ 
\frac{\|Bu\|}{\|u\|}\ \ge\ a^2, \\[0.5\baselineskip]
\beta_2 & := & 
\frac{\langle Vu',u'\rangle +a^2\langle Vu,u\rangle }{\langle Bu,u\rangle 
} = V(\xi)\frac{\langle u',u'\rangle +a^2\langle u,u\rangle }{\langle 
Bu,u\rangle } \nonumber\\[-0.5\baselineskip] 
\label{17b} \\[-0.5\baselineskip] 
& = & V(\xi)\in[V_{\min},V_{\max}]\quad\mbox{(for 
some $\xi\in(0,\infty)$)}, \nonumber \\[0.5\baselineskip]
\label{18}
\beta_3 & := & \frac{\langle V'u,u'\rangle }{\langle Bu,u\rangle },\quad 
|\beta_3|\le |V'|_{\max}
\frac{\|u\|\,\|u'\|}{\|u'\|^2+a^2\|u\|^2}\le \frac{|V'|_{\max}}{2a}.
\end{eqnarray}
From (\ref{15}) to (\ref{18}) we obtain the following eigenvalue
enclosure result, where sums and products of sets are to be understood
in the canonical sense.
\begin{theorem}
\label{encl1}
All eigenvalues of problem $($\ref{5}$)$ are contained in the set
$$
[a^2,\infty)+iaR[V_{\min},V_{\max}]+\frac{R}{2}|V'|_{\max}\cdot\Delta
,
$$
with $\Delta$ denoting the closed unit disc in $\C$.
\end{theorem}

Of course, separate bounds for the real and imaginary parts of the
eigenvalues can easily be extracted (with loss of information!) from
Theorem 3.1:
\begin{corollary}
\label{encl1a}
For each eigenvalue $\lambda$ of problem $($\ref{5}$)$,
\begin{eqnarray*}
a^2-\frac{R}{2}|V'|_{\max}\!\! & \le \!\! & \mbox{\rm Re}\,\lambda \\
aRV_{\min} -\frac{R}{2}|V'|_{\max}\!\! & \le \!\! & \mbox{\rm Im}\,\lambda\ 
\le\ aRV_{\max}+\frac{R}{2}|V'|_{\max} .
\end{eqnarray*}
\end{corollary}

\begin{figure}[h]
%\begin{center}  
  \epsfverbosetrue \hspace{1.0cm} \epsfxsize=10cm \epsfbox{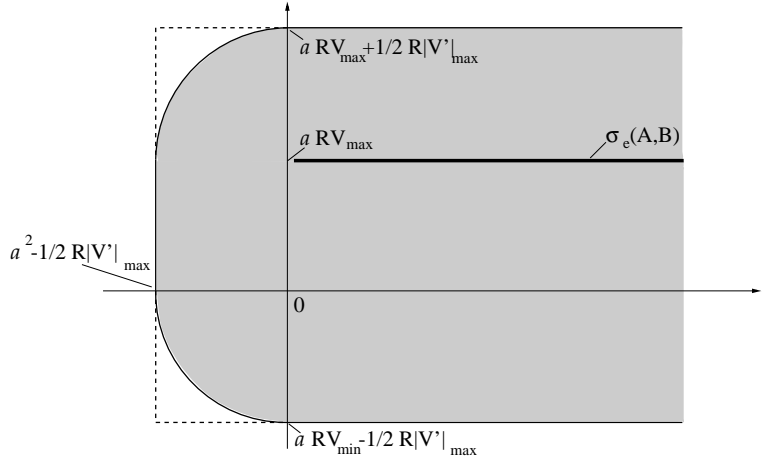}
%\end{center}
\caption{Illustration of Theorem 3.1 and Corollary 3.2} % !!!!!!!!!!!!!!!!!!   
\label{bild1}
\end{figure}

\begin{figure}[h]
%\begin{center}   
  \epsfverbosetrue \hspace{1.0cm} \epsfxsize=10cm \epsfbox{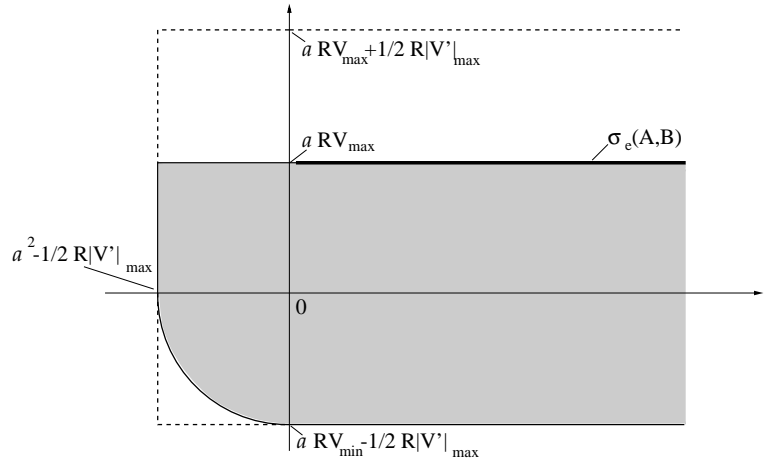}
%\end{center}
\caption{Illustration of Theorem 3.3}
\label{bild2}
\end{figure}

Figures \ref{bild1} and \ref{bild2} illustrate the results of Theorem
\ref{encl1}, Corollary \ref{encl1a} and Theorem \ref{encl2}.
Concretely, we have chosen $V$ to be the Blasius profile and $a=0.179,
R=580$ here.

The bounds given in Corollary 3.2 have been obtained in
\cite{fischer1} already, by more involved analytical means. These
results can be substantially improved by the following simple
calculation, if sign restrictions such as (\ref{4}) are used:
\begin{eqnarray*}
\mbox{\rm Re}\,\langle V'u,u'\rangle & = & 
\frac{1}{2}(\langle V'u,u'\rangle +\overline{\langle V'u,u'\rangle })\ =\
\frac{1}{2}\int_0^\infty V'(|u|^2)'\, dx\\
& = & -\frac{1}{2}\int_0^\infty V''|u|^2\, dx\ =\ 
-\frac{1}{2}V''(\xi)\langle u,u\rangle  \quad\mbox{(for some 
$\xi\in(0,\infty)$)} 
\end{eqnarray*}
so that $\beta_3$ defined in (\ref{18}) satisfies
\begin{equation} \label{19}
\mbox{\rm Re}\,\beta_3=-\frac{1}{2}V''(\xi)\,\frac{\langle u,u\rangle 
}{\langle Bu,u\rangle } \ \in\ 
-\frac{1}{2}[V''_{\min},V''_{\max}]\cdot[0,\frac{1}{a^2}] \end{equation}
which possibly restricts the enclosure set for $\beta_3$ used in Theorem 3.1
(which was obtained from (\ref{18})). Instead of formulating a
corresponding improved theorem for the most general situation, we
concentrate on the case where $V''\le 0$, which is true e.g. for the
Blasius profile (see (\ref{4})). Then (\ref{19}) yields Re$\,\beta_3\ge 0$,
which together with (\ref{15}) to (\ref{18}) provides the following
improvement of Theorem 3.1:
\begin{theorem}
\label{encl2}
If $V''\le 0$, all eigenvalues of problem $($\ref{5}$)$ are contained
in the set
$$
[a^2,\infty)+iaR[V_{\min},V_{\max}]+\frac{R}{2}|V'|_{\max}\cdot\Delta^-
,
$$
where $\Delta^-:=\{z\in\Delta\, :\, \mbox{\rm Im}\,z\le 0\}.$
\end{theorem}

As in Corollary 3.2, we can extract separate bounds for real and
imaginary parts of the eigenvalues, which provides here the improved
bound
$$
aRV_{\min} -\frac{R}{2}|V'|_{\max}\ \le\ \mbox{\rm Im}\,\lambda\ 
\le\ aRV_{\max} .
$$
In particular all eigenvalues lie "below" the essential spectrum,
i.e.  have imaginary parts less or equal $aR$.

%\clearpage

\begin{acknowledgements}
The authors are grateful to Andrei Shkalikov for helpful comments and 
suggestions.
\end{acknowledgements}

\address{Mathematisches Institut I, University of Karlsruhe, D-76128
  Karlsruhe}

\end{document}